\documentclass[12pt,a4paper,leqno,twoside]{article}
\usepackage{amsfonts,latexsym,bbm}

\usepackage{amsmath}
\usepackage{amssymb}
\usepackage{amscd}
\usepackage{amsopn}
\usepackage{eurosym}
\usepackage{amsthm}

\makeindex
\makeatletter
\@addtoreset{equation}{section}
\makeatother

\newtheorem{lemma}{Lemma}[section]
\newtheorem{proposition}[lemma]{Proposition}
\newtheorem{corollary}[lemma]{Corollary}
\newtheorem{theorem}[lemma]{Theorem}
\newtheorem{para}[lemma]{}
\newtheorem{remark}[lemma]{Remark}

\newtheorem{definition}[lemma]{Definition}

\newtheorem{counterexample}[lemma]{Counterexample}

\newcommand{\beweis}{\begin{proof}[Beweis]}
\newcommand{\beweisende}{\end{proof}}

\newcommand{\real}{\mathbbm{R}}

\begin{document}

\title{Proper actions and proper invariant metrics}

\author{H.~Abels, A.~Manoussos and G.~Noskov
\thanks{During this research the second and the third named authors were
  fully supported by the SFB 701 ``Spektrale Strukturen und
  Topologische
   Methoden in der Mathematik" at the
 University of Bielefeld, Germany. They are grateful for its
 generosity and hospitality. The paper was finished while the first
 named author was staying at MSRI in Berkeley. He wishes to thank the
 MSRI for its hospitality and support, as well as the SFB 701 for its
 support.}
}

%\author[H.~Abels]{Herbert Abels}
%\address[H.~Abels]{Fakult\"{a}t f\"{u}r Mathematik, Universit\"{a}t
% Bielefeld, Postfach 100131, D-33501 Bielefeld, Germany}
%\email[H.~Abels]{abels@math.uni-bielefeld.de}

% \author[A. Manoussos]{Antonios Manoussos}
% \address{Fakult\"{a}t f\"{u}r Mathematik, SFB 701, Universit\"{a}t
% Bielefeld, Postfach 100131, D-33501 Bielefeld, Germany}
% \email{amanouss@math.uni-bielefeld.de}

% \author[G.A. Noskov]{Gennady A. Noskov}
% \ address{Sobolev Institute of Mathematics, Pevtsova 13, Omsk 644099,
% Russia and%Fakult\"{a}t f\"{u}r Mathematik, Universit\"{a}t
% Bielefeld, Postfach 100131, D-33501 Bielefeld, Germany}
% \email{noskov@math.uni-bielefeld.de}

\date{}

\maketitle

\begin{abstract}
We show that if a (locally compact) group $G$ acts properly on a
locally compact $\sigma$-compact space $X$ then there is a family of
$G$-invariant proper continuous finite-valued pseudometrics which
induces the topology of $X$. If $X$ is furthermore metrizable then $G$
acts properly on $X$ if and only if there exists a $G$-invariant
proper compatible metric on $X$.
\end{abstract}

\vglue0.3cm {\scshape Subject classification} [2000]: Primary 37B05, 54H15; Secondary 54H20, 54D45.
\vglue0.3cm {\scshape Keywords:} Proper action, group of
isometries, proper metric, proper pseudometric, Heine-Borel metric.

\section{Introduction}\label{sec1}

We establish a close connection between proper group actions and
groups of isometries. There is an old result in this direction, proved
in 1928 by van Dantzig and van
der Waerden %\cite{d-w}.
It says that for a locally compact connected
metric space $(X,d)$ its group $G=\mbox{Iso}(X,d)$ of isometries is locally
compact and acts properly.
That the action is proper is no longer true in general, if $X$ is not
connected, although $G$ is sometimes still locally compact, see
\cite{m-s}. Concerning properness of the action, Gao and Kechris \cite{kechris}
proved the following result. If $(X,d)$ is a proper metric space, then
$G$ (is locally compact and) acts properly on $X$. Recall that a
metric $d$ on a space $X$ is called proper if balls of bounded radius
have compact closures.

There is the following converse result. If a locally compact group $G$
acts properly on a locally compact $\sigma$-compact metrizable space
$X$, then there is a compatible $G$-invariant metric $d$ on $X$
\cite{koszul}. In this paper we prove that under these hypotheses
there is actually a compatible $G$-invariant {\em proper} metric on
$X$. We call a metric on a topological space compatible if induces its
topology. Note that a proper metric space is $\sigma$-compact. For
the records, here is one version of our main result, namely the one
for metrizable spaces (see also theorem \ref{theo4.2}).

\begin{theorem}\label{theo1.1}
Suppose the (locally compact) topological group $G$ acts properly on
the metrizable locally compact $\sigma$-compact topological space
$X$. Then there is a
$G$-invariant proper compatible metric on $X$.
\end{theorem}

These results raise the question if they generalize to the
non-metri\-zable case. We give a complete answer as follows. Recall
that a pseudometric on $X$ is a
function $d$ on $X\times X$ which has all the properties of a metric,
except that its value may be $\infty$ and that $d(x,y)=0$ may not
imply that $x=y$. For a precise definition see below definition
\ref{def2.1}.
A locally compact space is $\sigma$-compact if and only if
has a proper finite-valued continuous pseudometric, as is
easily seen, see e.g.~below, the proof of corollary \ref{cor5.3}.
It then actually has a family of such pseudometrics
which induces the topology of $X$. The corresponding statement for the
equivariant situation is the following version of the main result of
our paper, namely for not necessarily metrizable spaces (see also
theorem \ref{theo4.1}).

\begin{theorem}\label{theo1.2}
Let $G$ be a (locally compact) topological group which acts properly
on a locally compact $\sigma$-compact Hausdorff space $X$. Then there
is a family of $G$-invariant proper finite-valued continuous
pseudometrics on $X$ which induces the topology of $X$.
\end{theorem}

The connection of theorem  \ref{theo1.1} and theorem \ref{theo1.2} is
given by the following result. We are in the case of theorem
\ref{theo1.1} iff there is a countable family as in theorem
\ref{theo1.2}. For a precise statement see corollary \ref{cor4.4}.

Note that continuity of the pseudometrics follows from the other
properties, see remark \ref{rem5.5}.

This theorem may be considered as the converse of the following
theorem, see below theorem \ref{theo3.1}.

\begin{theorem}\label{theo1.3}
Let $X$ be a topological space and let $\mathcal{D}$ be a family of
proper continuous finite-valued pseudometrics on $X$, which induces
the topology of $X$. Let $G$ be the
group of all bijective maps $X\to X$, leaving every $d\in\mathcal D$
invariant. Endow $G$ with the compact--open topology. Then $G$ is a
locally compact topological group and acts properly on $X$.
\end{theorem}

The main result of our paper has been proved already for the special case of a smooth manifold. Namely Kankaanrinta proved in \cite{kaa1} that if a Lie group $G$ acts
properly and smoothly on a smooth manifold $M$, then $M$ admits a complete $G$-invariant Riemannian metric. A consequence of our main result for the metrizable case
is the following result of Haagerup and Przybyszewska \cite{haag}: Every second countable locally compact group has a left invariant compatible proper metric which
generates its topology, see below corollary \ref{cor9.5}. Proper $G$--invariant metrics have been used in several fields of mathematics, see \cite{jan} and
\cite{skand}. For more information about related work, open questions and miscellaneous remarks see the last chapter of this paper.

\section{Preliminaries}\label{sec2}
\noindent
{\large\bf Pseudometrics}
\vglue0.4cm
\begin{definition}\label{def2.1}
A pseudometric $d$ on a set $X$ is a function $d:X\times
X\to[0,+\infty]$ which fulfills for $x,y,z\in X$ the following
properties
\begin{itemize}
\item[i)] $d(x,x)=0$,
\item[ii)] $d(x,y)=d(y,x)$,
\item[iii)] $d(x,y)+d(y,z)\ge d(x,z)$.
\end{itemize}
\end{definition}
Thus, loosely speaking, a pseudometric is a metric except that its
values may be $+\infty$ and $d(x,y)=0$ does not imply $x=y$. A family
$\mathcal D$ of pseudometrics on $X$ induces a topology on $X$, for
which finite intersections of balls $B_d(x,r):=\{y\in X; d(x,y)<r\}$
with $x\in X$, $d\in\mathcal D$ and $r\in [0,\infty)$ form a
basis. This topology is the coarsest topology for which every
$d\in\mathcal D$ is a continuous function on $X\times X$. The topology
of a topological space $X$ is induced by a family of pseudometrics if
and only if $X$ is completely regular, see
\cite[Ch. X, \S 1.4 Theorem 1 and \S 1.5 Theorem 2]{bour}. A topological
space $X$ is called
metrizable if its topology is
induced by an appropriately chosen metric $d$ on $X$. Such a metric
$d$ on $X$ is then called a {\em compatible metric}.

{\em From now on we will call a locally compact Hausdorff space simply a ``space'', for short}. Recall that a space is called $\sigma$--compact if it can be written
as a countable union of compact subsets. A $\sigma$--compact space is metrizable if and only if it is second countable, i.e., its topology has a countable base, see
\cite[Ch. IX, \S 2.9 Corollary]{bour}.

A pseudometric $d$ on a space $X$ will be called proper if every ball
of finite radius has compact closure. A space $X$ together with a
compatible proper metric $d$ will be called a {\em proper metric
  space}. It is also called a {\em Heine--Borel space} by some
authors and also a {\em finitely--compact space} by others. Important
examples of proper metric spaces are the Euclidean
spaces and the space $\mathbb Q_p$ of rational $p$--adics with their
usual metrics.

The topology of a space can be induced by a family of pseudometrics,
since a space (understood: locally compact Hausdorff) is completely
regular. The topology of a $\sigma$--compact space can be induced by a
family of
proper finite--valued pseudometrics (see corollary \ref{cor5.3}).
One of our main results, theorem
\ref{theo1.2}, spells out for which actions there is a family of
invariant proper finite--valued pseudometrics inducing the topology,
namely the proper actions. And theorem \ref{theo1.3} says that these
are essentially the only ones for which such a family exists.

Now let $(X,\mathcal D)$ be a space $X$ together with a family
$\mathcal D$ of pseudometrics inducing its topology. A case of
particular importance is when $\mathcal D$ consists of just one
metric, which by assumption induces the topology of $X$. Let
$G=\mbox{Iso}(X,\mathcal D)$ be the group of isometries of
$(X,\mathcal D)$, that is the group of all bijections $X\to X$
leaving every $d\in\mathcal D$ invariant. Endow $G$ with the topology
of pointwise convergence. Then $G$ will be a topological group
\cite[Ch. X, \S 3.5 Corollary]{bour}. On $G$ there is also the topology of
uniform convergence on
compact subsets which is the same as the compact--open topology. In our
case, these topologies coincide with the topology of pointwise
convergence, and the natural action of $G$ on $X$ is continuous
\cite[Ch. X, \S 2.4 Theorem 1 and \S 3.4 Corollary 1]{bour}. We shall
prove soon, that if at least one of the
pseudometrics $d$ in $\mathcal D$ is proper then $G$ is locally
compact. In this case the natural action of $G$ on $X$ is even
proper. We will discuss this notion now.

\begin{definition}\label{def2.2}
A continuous map $f:X\to Y$ between spaces is called proper if one of
the following two equivalent conditions holds
\begin{itemize}
\item[i)] $f^{-1}(K)$ is compact for every compact subset $K$ of $Y$.
\item[ii)] $f$ is a closed map and the inverse image of every
  singleton is compact.
\end{itemize}
\end{definition}

Let $G$ be a topological group. Suppose a continuous action of $G$ on
a space $X$ is given.

\begin{proposition}\label{prop2.3}
{\bf and Definition}
The following conditions are equivalent
\begin{itemize}
\item[i)] The map $G\times X\longrightarrow X\times X$, $(g,x)\longmapsto
  (gx,x)$, is proper.
\item[ii)] For every pair $A$ and $B$ of compact subsets of $X$ the
  transporter
$$
G_{AB}:=\{g\in G;\ g A\cap B\ne \emptyset\}
$$
from $A$ to $B$ is compact.
\item[iii)] Whenever we have two nets $(g_i)_{i\in I}$ in $G$ and
  $(x_i)_{i\in I}$ in $X$, for which both $(x_i)_{i\in I}$ and
  $(g_ix_i)_{i\in I}$ converge, then the net $(g_i)_{i\in I}$ has a
  convergent subnet.
\end{itemize}
The action of $G$ on $X$ is called proper if one of these conditions holds.
\end{proposition}

For a proof see \cite[Ch. I, \S 10.2 Theorem 1 and Ch. III, \S 4.4
Proposition 7]{bour}. For more information on proper group actions see
the forthcoming book \cite{AS}.

Note that if the action of $G$ on $X$ is proper, then $G$ is locally
compact, by ii). And if furthermore $X$ is $\sigma$--compact, then
$G$ is also $\sigma$--compact, by ii).

It is useful to rephrase the definition of properness in terms of
limit sets. Let $(x_i)_{i\in I}$ be a net in the -- not necessarily
locally compact -- topological space $X$. We say that the net
$(x_i)_{i\in I}$ diverges and write $x_i \underset{i\in
  I}\longrightarrow \infty$, if the
net $(x_i)_{i\in I}$ has no convergent subnet. If $X$ is locally
compact, a net $(x_i)_{i\in I}$ in $X$ diverges if and only if it
converges to the additional point $\infty$ of the one point (also called
Alexandrov--) compactification of $X$.

Let again the topological group $G$ act on the space $X$. For $x\in X$
the {\em limit set} $L(x)$ is defined by
\begin{align*}
L(x):=&\{y;\text{\em there exists a divergent net }(g_i)_{i\in I}\
\text{\em in } G \\
&\text{\em such that } (g_i x)_{i\in I}\ \text{\em converges to }y\}
\end{align*}
and the {\em extended limit set} $J(x)$ is defined by
\begin{align*}
J(x):=&\{y;\text{\em there exists a divergent net }(g_i)_{i\in I}
\ \text{\em in } G \\
&\text{\em  and a net }(x_i)_{i\in I}\ \text{\em  in } X
\text{\em converging to } x, \\
&\text{\em such that }(g_i x_i)_{i\in I}\ \text{\em  converges to }y\}.
\end{align*}

Thus, the action of $G$ on $X$ is proper if and only if the following
condition holds:
\begin{itemize}
\item[iv)] $J(x)=\emptyset$ {\em for every} $x\in X$,\\
since iv) is equivalent to iii). If furthermore
$\mathcal D$ is a family of pseudometrics inducing the topology of
$X$ and every $g\in G$ leaves every $d\in\mathcal D$ invariant, then
it is easy to see that
\item[v)] $L(x)=\emptyset$ {\em implies} $J(x)=\emptyset$.
\end{itemize}

\section{The group of isometries of a proper metric space}\label{sec3}

Let again $X$ be a locally compact Hausdorff space, let $\mathcal D$
be a family of pseudometrics inducing the topology of $X$ and let $G$
be the group of isometries of $(X,\mathcal D)$ with its natural
topology, as above.

\begin{theorem}\label{theo3.1}
If at least one of the pseudometrics in $\mathcal D$ is proper then
$G$ is locally compact and the natural action of $G$ on $X$ is proper.
\end{theorem}

The special case that $\mathcal D$ consists of just one metric is due
to Gao and Kechris \cite{kechris}, as follows.

\begin{theorem}\label{theo3.2}
If $(X,d)$ is a proper metric space then its group $G$ of isometries
is locally compact and its natural action of $G$ on $X$ is proper.
\end{theorem}

\noindent
{\bf Proof of theorem \ref{theo3.1}}. It suffices to show that the
natural action of $G$ on $X$ is proper. To prove this we will show
that the limit set $L(x)$ is empty for every $x\in X$. Thus let
$(g_i)_{i\in I}$ be a net in $G$ for which $(g_i x)_{i\in I}$ converges
to a point, say $y$, in $X$. We have to show that the net $(g_i)_{i\in
  I}$ has a convergent subnet. We may assume that $g_i x$ is contained
in the relatively compact ball $B_d(y,r)$ for every $i\in I$, where
$d$ is a proper pseudometric in $\mathcal D$ and $r>0$. We will use the
Arzela--Ascoli theorem. Let $z\in X$. The points $g_i z$, $i\in I$,
are contained in the ball $B_d(z,R)$, where $R=r+d(x,z)$. Thus the set
$\{g_i z; i\in I\}$ is relatively compact for every $z\in X$. The
family of maps $\{g_i; i\in I\}$ is uniformly equicontinuous, being a
subset of the uniformly equicontinuous family $G$ of maps from $X$ to
$X$. It follows from the Arzela--Ascoli theorem that the net
$(g_i)_{i\in I}$ has a subnet $(g_j)_{j\in J}$ which converges
uniformly on compact subsets to a map $g$. Clearly, $g$ leaves every
$d\in\mathcal D$ invariant. To see that $g$ is actually invertible
look at the net $(g_j^{-1})_{j\in J}$. We have $g_j^{-1}y\in B_d(x,r)$
and hence $g_j^{-1}z\in B_d(z,R')$ where $R'=r+d(x,z)$. Then the net
$(g_j^{-1})_{j\in J}$ has a subnet which converges uniformly on
compact subsets to a map $f$. It then follows that $f$ and $g$ are
inverse of each other.

\begin{remark}\label{rem3.3}
{\em The sets $K(E):=\{ x\in X\,;\,\, Ex\,\,\mbox{is relatively
  compact}\}$, where \linebreak
$E\subset \mbox{Iso}(X,d)$ played a crucial role in
\cite{m-s} where it is proved that they are open--closed subsets of
$X$. In the case of a proper metric space $(X,d)$ the set $K(E)$ is
either the empty set or the whole space $X$ as shown in the proof of
Theorem \ref{theo3.1}. Using Bourbaki \cite[Ch. X, Exercise 13, p.~323]{bour}
we may also show that sets $K(E)$ are open-closed
subsets of $X$ but we must be careful! Even
in the legendary ``Topologie G\'{e}n\'{e}rale'' of Bourbaki there is at
least one mistake! Precisely in the aforementioned Exercise 13 of
Ch. X, p. 323, part d) it is said that if $E$ is a uniformly
equicontinuous family of homeomorphisms of a locally compact uniform
space $X$ then $K(E)$ is a closed subset of $X$. This is not true if
$E$ is not a subset of a uniformly equicontinuous \textit{group} of
homeomorphisms of $X$ as we can easily see by the following
counterexample.}
\end{remark}

\begin{counterexample}\label{cexam3.4}
{\rm Let
\[
X=\bigcup_{k=1}^{\infty}\{ (x,y)\,;\, x=\frac{1}{k},\,y\geq 0\}\cup \{
(x,y)\,;\, x=0,\, y>0\}
\]
be endowed with the Euclidean metric. Consider the family
$E=\{f_{n}\}$ of selfmaps of $X$ defined by
$f_{n}(x,y)=(x,\frac{y}{n})$. The family $E$ consists of
uniformly equicontinuous homeomorphisms of $X$ and
$K(E)=\bigcup_{k=1}^{\infty}\{ (x,y)\,;\, x=\frac{1}{k},\,y\geq 0\}$
as can be easily checked. Hence the set
$K(E)$ is not closed in $X$.}
\end{counterexample}

\section{Proper invariant metrics and pseudometrics, outline of the
  proof}\label{sec4}

The main results of our paper are the following converses of theorems
\ref{theo3.1} and \ref{theo3.2}. Again, $X$ is a space, i.e., a
locally compact Hausdorff space and $G$ is a Hausdorff topological
group. Suppose we are given a continuous action of $G$ on $X$.

\begin{theorem}\label{theo4.1}
Suppose $X$ is $\sigma$--compact. If the action of $G$ on $X$ is
proper then there is a family
$\mathcal D$ of proper finite--valued $G$--invariant pseudometrics on
$X$, which induces the topology of $X$.
\end{theorem}

\begin{theorem}\label{theo4.2}
Suppose $X$ is $\sigma$--compact. If the action of $G$ on $X$ is
proper and $X$ is metrizable, then
there is a compatible $G$--invariant proper metric $d$ on $X$.
\end{theorem}

\begin{remark}\label{rem4.3}
If the action is proper, it is easy to see that the kernel of the
action $K:=\{g\in G\ ;\ gx=x\text{ for every }x\in X\}$ is compact and
the action map induces an isomorphism of topological groups of $G/K$
onto a closed subgroup of $\mbox{Iso}(X,\mathcal D)$,
resp. $\mbox{Iso}(X,d)$. We thus have a complete correspondence
between proper actions and isometry groups of proper metrics or pseudometrics.
\end{remark}

\begin{corollary}\label{cor4.4}
Suppose $X$ is $\sigma$--compact and $G$ acts properly on $X$. Then
the following properties of $X$ are equivalent
\begin{itemize}
\item[a)] $X$ is metrizable.
\item[b)] There is a compatible $G$--invariant proper metric on $X$.
\item[c)] There is a countable family of finite--valued
  pseudo--metrics on $X$, which induces the topology of $X$.
\item[d)] There is a countable family of proper finite--valued
  $G$--invariant pseudometrics on $X$, which induces the topology of $X$.
\end{itemize}
\end{corollary}

\begin{proof}
a) $\Longrightarrow$ b) by theorem \ref{theo4.2}, b) $\Longrightarrow$
d) and d) $\Longrightarrow$ c) are trivial, c) $\Longrightarrow$ a) is
a well known theorem of topology \cite[Ch. IX, \S 2.4 Corollary 1]{bour}
whose proof is similar to the argument in the last paragraph of the
proof of lemma \ref{lem8.10} a).
\end{proof}

The proof of theorems \ref{theo4.1} and \ref{theo4.2} will occupy most
of the remainder of the
paper. Let us briefly describe the plan of the proof. We describe the
plan for the case of a family of pseudometrics, the proof for the
metrizable case simplifies at some points.
\begin{enumerate}
\item We first construct a family $\mathcal D$ of pseudometrics on
  $X$, with values in [0,1] which induces the topology of $X$, see
  section \ref{sec5}.
\item Next we show how to make every $d\in\mathcal D$ $G$--invariant,
  see section \ref{sec6}.
\item Then we make every $d\in\mathcal D$ orbitwise proper, see
  section \ref{sec7}.
\item These steps are fairly routine. We then present our main tool,
  namely the ``measuring stick construction''.
Imagine a family of measuring sticks given by distances of closely
neighboring points. We then define a pseudometric $d$ on $X$ by taking
for $x,y$ in $X$ as $d(x,y)$ the infimum of all measurements along
sequences of points $x=x_0,\dots,x_n=y$ such that the distance of any
two adjacent points is given by measuring sticks. For a precise
definition, actually several equivalent ones, see section \ref{sec8}.
It turns out that we then get for an appropriate family of measuring
sticks a proper pseudometric. The disadvantage of this construction is
that there may be points which cannot be connected by sequences as
above. Equivalently, there may be points $x,y$ with $d(x,y)=\infty$.
\item We then use our ``bridge construction'', see section
  \ref{sec9}. Think of pairs of
  points with $d(x,y)<\infty$ as lying on the same island.
Thus what we call an island is an equivalence class of the equivalence
relation defined as $x\sim y$ iff $d(x,y)<\infty$.
 We connect (some of) these
  islands by bridges and attribute (high) weights to these bridges. We
  then define a new pseudometric similarly as above using the already
  defined pseudometric on the islands and the weights of bridges. We
  thus obtain a proper pseudometric with finite values and actually a
  whole family of such, which induces the topology of $X$. All these
  constructions are done in a $G$--invariant way, so that the
  resulting pseudometrics are $G$--invariant.
\end{enumerate}

\section{A compatible metric and proper pseudometrics}\label{sec5}

Again, by a space we mean a locally compact Hausdorff space. Recall
the following basic metrization result, see
\cite[Ch. IX, \S 2.9 Corollary]{bour}.

\begin{theorem}\label{theo5.1}
For a space $X$ the following properties are equivalent
\begin{itemize}
\item[a)] $X$ is second countable, i.e., its topology has a countable
  base.
\item[b)] The one--point compactification $\overline X$ of $X$ is
  metrizable.
\item[c)] $X$ is metrizable and $\sigma$--compact.
\end{itemize}
\end{theorem}
If a space is metrizable we may assume that the metric $d$ inducing
the topology has
values in [0,1]. We just have to replace $d$ by $d_1$ with
$d_1(x,y):=\frac{d(x,y)}{1+d(x,y)}$.

For the general case of a not necessarily metrizable $\sigma$--compact
space --- and for later use --- we need the following easy lemma,
whose proof is left to the reader.

\begin{lemma}\label{lem5.2}
A space $X$ is $\sigma$--compact if and only if there is a proper
continuous function $f:X\longrightarrow [0,\infty)$.
\end{lemma}

\begin{corollary}\label{cor5.3}
On every $\sigma$--compact space $X$ there is a family $\mathcal D$
of proper finite--valued pseudometrics inducing the topology of $X$.
\end{corollary}

\begin{proof}
Let $\mathcal D_0$ be the family of pseudometrics on $X$ of the form
$$
d_f(x,y):=|f(x)-f(y)|
$$
for $x,y\in X$, where $f:X\longrightarrow \mathbb R$ is a continuous
function. Then $\mathcal D_0$ induces the topology of $X$. Here we do
not use that $X$ is $\sigma$--compact. But in the next step we do. If
$X$ is $\sigma$--compact let $\mathcal D$ be the family $\mathcal
D:=\{d+d_f\ ;\ d\in\mathcal D_0\}$, where $f:X\longrightarrow\mathbb
R$ is proper and continuous. Then $\mathcal D$ induces the topology of
$X$ and consists of proper finite--valued pseudometrics.
\end{proof}
The same trick yields the following corollary.

\begin{corollary}\label{cor5.4}
The following properties of a space $X$ are equivalent.
\begin{itemize}
\item[a)] $X$ has a compatible proper metric.
\item[b)] $X$ is metrizable and $\sigma$--compact.
\item[c)] $X$ is metrizable and separable.
\item[d)] $X$ is second countable.
\end{itemize}
\end{corollary}

\begin{remark}\label{rem5.5}
{\rm Note the if a pseudometric $d$ belongs to a family of pseudometrics
inducing the topology of $X$ then $d$ is continuous. Since then
$B_d(x,r)$ is a neighborhood of $x$ for every $x\in X$ and every
$r>0$, and hence the function $y\longmapsto d(x,y)$ is continuous at
$x$ for every $x\in X$, which easily implies that $d$ is continuous by
the triangle inequality.}
\end{remark}

\section{Making the metrics or pseudometrics
  $G$--invariant}\label{sec6}

Now suppose $X$ is a space, $G$ is a Hausdorff topological group and a
proper continuous action of $G$ on $X$ is given.
\vglue0.3cm
\noindent
{\bf Step 2.}\ \ {\em If $X$ is $\sigma$--compact, then there is a family of
$G$--invariant continuous finite--valued pseudometrics inducing the
topology of $X$. If $X$ is furthermore metrizable then there is a
compatible $G$--invariant metric on $X$.}

We present two proofs.

The first one is due to Koszul \cite{koszul} and uses the concept of a
fundamental set, a concept we will need again, later on. The second
one uses the notion of an equicontinuous
action on the one--point compactification of $X$. Unfortunately, in
the process we loose the property that our (pseudo--)metrics are proper.

\begin{definition}\label{def6.1}
A subset $F$ of $X$ is called a {\em fundamental set} for the action of
$G$ on $X$ if the following two conditions hold.
\begin{itemize}
\item[a)] $GF=X$\
\item[b)] $G_{KF}$ has compact closure for every compact subset $K$ of $X$.
\end{itemize}
\end{definition}
Concerning b), recall the definition of the transporter $G_{AB}=\{g\in
G\ ; \ g\,A\cap B\ne \emptyset\}$ from $A$ to $B$. Note that only proper
actions can have a fundamental set, since a) implies that
$$
G_{AB}\subset G^{-1}_{BF} \cdot G_{AF}
$$
and hence $G_{AB}$ is relatively compact if $A$ and $B$ are compact,
by b), and then $G_{AB}$ is actually compact, by continuity of the
action. There is the following converse, see \cite{koszul}.

\begin{proposition}\label{prop5.2}
If $X$ is $\sigma$--compact, then there is an open fundamental set for
every proper action.
\end{proposition}
\vglue0.3cm
\noindent
{\bf Step 2, 1$^{\text{\bf st}}$ proof.}\ \ Let $F$ be an open fundamental set
for the action of $G$ on $X$. Let $d$ be a continuous finite--valued
pseudometric on $X$. Let $d'$ be the supremum of all pseudometrics on
$X$ with the property that $d'\mid F\times F\le d$ and
$d'\mid(X\smallsetminus F)\times (X\smallsetminus F)=0$. Explicitly, let $r$
be the function on $X$ with $r_d(x)=d(x,X\smallsetminus F):=\inf \{d(x,y)\
;\ y\in X\smallsetminus F\}$. Then
$$
d'(x,y)=\min\{d(x,y)\ ,\ r_d(x)+r_d(y)\}.
$$
Note that for every $x\in F$ there is a neighborhood of $x$ where $d$
and $d'$ coincide. The function $d'$ is a finite--valued continuous
pseudometric and the
function $G\longrightarrow \mathbb R$, $g\longmapsto d'(gx,gy)$ is
continuous and has compact support, namely contained in
$G_{\{x,y\},F}$. Define
$$
d''(x,y)=\int_G d'(gx,gy) dg
$$
where $dg$ is a right invariant Haar measure on $G$. Then $d''$ is a
$G$--invariant pseudometric on $X$. The pseudometric $d''$ is actually
a metric if $d$ is a metric. Furthermore $d''$ is continuous for every
$d\in\mathcal D$, by a uniform equicontinuity argument for functions
on compact spaces. Thus the family $\mathcal D'' =\{d''\ ;\
d\in\mathcal D\}$ induces a weaker topology than $\mathcal D$. The
two topologies are actually equal since for every neighborhood $V$ of
$x\in X$ there are a compact neighborhood $V_1$ of $x$ in $X$ and a
compact neighborhood $U_1$ of $e$ in $G$ such that $U_1 V_1\subset V$
and $U_1(X\smallsetminus V)\subset X\smallsetminus V_1$ and hence
$$
d''(x,y)\ge d'(x,X\smallsetminus V_1)\cdot\int_{U_1}dg
$$
for every $y\in X\smallsetminus V$, which implies our claim for $x\in
F$ and hence for every $x$ by $G$--invariance of the two topologies.
\vglue0.3cm
\noindent
{\bf ${\mathbf 2}^{\text{\bf nd}}$ proof.}\ \ This proof is based on the notion
of an equicontinuous group action. Consider the
one point compactification $\overline X=X\cup\{\infty\}$. The action
of $G$ on $X$ extends to an action of $G$ on $\overline X$ by defining
$g(\infty)=\infty$ for every $g\in G$. The extended action is
continuous. Let $\mathcal D$ be a family of pseudometrics on
$\overline X$ which induces the topology of $\overline X$. Without
further assumptions on $X$ we can take the family $\{d_f; f:\overline
X\to [0,1] \mbox{ continuous}\}$, see the proof of corollary
\ref{cor5.3}. If $\overline
X$ is metrizable, we can take $\mathcal D$ to consist of just one
element. This is the case if and only if $X$ is metrizable and
$\sigma$--compact, see theorem \ref{theo5.1}. In any case, define for
$d\in\mathcal D$ and $x,y\in X$
$$
d'(x,y):=\sup_{g\in G} d(gx,gy),
$$
and set $\mathcal D'=\{d'\ ;\ d\in\mathcal D\}$. We claim that
$\mathcal D'$ induces the topology of $X$. Obviously, the topology
induced by $\mathcal D'$ is finer than the topology of $X$, since
$d'\ge d$ and $\mathcal D$ induces the topology of $X$.

Concerning the converse, consider the following property. The action
of $G$ on $X$ is called {\em pointwise equicontinuous with respect} to
$\mathcal D$ if for every $x\in X$, $d\in\mathcal D$ and $\epsilon
>0$ there is a neighborhood $U$ of $x$ such that for $y\in U$ we have
$d(gx,gy)<\epsilon$ for every $g\in G$. Clearly, if this holds the
topology defined by $\mathcal D'$ is weaker than the topology of $X$
and our claim is proved. It thus remains to show

\begin{lemma}\label{lem6.3}
Let $X$ be a space and let $G$ be a topological group acting properly
on $X$. Let $\mathcal D$ be a family of pseudometrics on $\overline
X$ inducing the topology of $\overline X$. Then $G$ acts pointwise
equicontinuously on $X$ with respect to $\mathcal D$.
\end{lemma}

\begin{proof}
Arguing by contradiction, assume that there are $d\in\mathcal D$,
$x\in X$, $\epsilon >0$ and a net $(x_i)_{i\in I}$ in $X$ converging
to $x$ and a net $(g_i)_{i\in I}$ in $G$ such that $d(g_i x,
g_i x_i)\ge \epsilon$ for every $i\in I$. It follows that
$g_i\longrightarrow \infty$, since otherwise the net $(g_i)_{i\in I}$
has a convergent subnet, say $(g_j)_{j\in J}$ converging to $g\in
G$. Then $g_jx\underset{j\in J}\longrightarrow gx$ and $g_j
x_j\underset{j\in J}\longrightarrow gx$ contradicting $d(g_i x, g_i
x_i)\ge\epsilon$ for every $i\in I$. It follows next that $g_i
x_i\underset{i\in I}\longrightarrow\infty$, since otherwise there
would be a subnet $(g_jx_j)_{j\in J}$ converging to a point of $X$,
  which implies that there would be a convergent subnet of
  $(g_j)_{j\in J}$, by properness of the action. Thus $g_i
    x_i\underset{i\in I}\longrightarrow\infty$ and $g_i\underset{i\in
      I}\longrightarrow\infty$, which implies $g_i x\underset{i\in
      I}\longrightarrow\infty$, again by properness of the action. But
    then $d(g_i x, g_ix_i)\underset{i\in I}\longrightarrow~0$, since
    $d$ is continuous on $\overline X$. This contradicts our
    assumption and finishes the proof.
\end{proof}

\begin{remark}\label{rem6.4}
{\em The $2^{\text{nd}}$ proof shows step 2 for the metrizable case only
  under the additional assumption that $\overline X$ is metrizable,
  i.e., that $X$ is metrizable and $\sigma$--compact. This is enough
  for our main results, though, because there all spaces are
  $\sigma$--compact.}
\end{remark}

\begin{remark}\label{rem6.5}
{\em The pseudometrics we obtain by these proofs are not proper, in
general. This is clear for the second proof. For the first proof, even
if we start from a proper (pseudo--) metric $d$, we obtain in case
that the orbit space $G\backslash X$ is compact -- so $F$ is
relatively compact -- that $d''$ has an upper bound.}
\end{remark}

\begin{remark}\label{rem6.6}
{\em One could rephrase the notion of pointwise equicontinuity in terms of
the unique uniformity on the compact space $\overline X$. We chose
here to use the language of pseudometrics since proper (pseudo--)
metrics are our final goal.}
\end{remark}

\section{Orbitwise proper metrics and pseudometrics}\label{sec7}

If $G$ acts on $X$ we denote by $\pi:X\longrightarrow G\backslash X$
the natural map to the orbit space. We will call a pseudometric $d$
on $X$ {\em orbitwise proper} if $\pi\left(B_d(x,r)\right)$ has compact closure
for every $x\in X$ and $0<r<\infty$. Again, we assume the notation and
hypotheses of the last section.
\vglue0.3cm
\noindent
{\bf Step 3.}\ \ {\em If $X$ is $\sigma$--compact there is a family of
$G$--invariant orbitwise proper finite--valued pseudometrics on $X$
inducing the topology of $X$. If $X$ is furthermore metrizable there
is a $G$--invariant orbitwise proper compatible metric on $X$.}

\begin{proof}
If $X$ is a space with a proper action, then the orbit space
$G\backslash X$ is Hausdorff as well, see
\cite{bour}. Clearly, $G\backslash X$ is locally compact. If
furthermore $X$ is $\sigma$--compact, so is $G\backslash X$. So
there is a proper continuous function $f:G\backslash X\to
[0,\infty)$, see lemma \ref{lem5.2}. The pseudometric
$d':=d_{f\circ\pi}$ on $X$ defined by
$$
d'(x,y)=\left| f\pi(x) - f\pi(y)\right|
$$
for $x,y\in X$ is orbitwise proper, continuous and
$G$--invariant. Hence if $\mathcal D$ is a family of finite--valued
$G$--invariant pseudometrics on $X$ inducing the topology of $X$, so
is $\mathcal D'=\{d+d'\ ;\ d\in\mathcal D\}$ and furthermore every
pseudometric of this family is orbitwise proper.
\end{proof}

\section{The measuring stick construction}\label{sec8}

We first present our measuring stick construction in three equivalent
ways. We then give a sufficient condition under which the resulting
pseudometric is proper. This will be applied to our situation and
yields step 4 of our proof.

\begin{para}\label{para8.1}
{\em Let $X$ be a set, let $d$ be a pseudometric on $X$ and let $\mathcal
U$ be a covering of $X$. We then define a new pseudometric
$d'=d'(d,\mathcal U)$ on $X$ depending on $d$ and $\mathcal U$ as
follows: $d'$ is the supremum of all pseudometrics $d''$ on $X$ with
the property that $d''|U\times U\le d|U\times U$ for every
$U\in\mathcal U$.}
\end{para}

\begin{para}\label{para8.2}
{\em We think of pairs $(x,y)$ of points lying in one $U\in\mathcal U$ as}
measuring sticks {\em or} sticks, {\em for short. A sequence $x=x_0$,
$x_1,\dots,x_n=y$ of points in $X$, such that any two consecutive
points form a stick, will be called a} stick path {\em from $x$ to $y$
of length $n$ and $d$--length $\sum^n_{i=1}d(x_{i-1},x_i)$. We claim that
$d'(x,y)$ is the infimum of $d$--lengths of all stick paths from $x$
to $y$. Since on one hand defining $d'$ in this way clearly gives a
pseudometric on $X$ and $d'|U\times U\le d|U\times U$.
And, on the other hand, for every pseudometric $d''$ with the two
properties above we have that $d''(x,y)$ is at most equal to the
$d$--length of any stick path from $x$ to $y$, because for every stick
path $x=x_0,x_1,\dots,x_n=y$ we have
$$
d''(x,y)\le\sum^n_{i=1}d''(x_{i-1},x_i)\le\sum^n_{i=1}d (x_{i-1},x_i).
$$
We thus
obtain the following properties of $d'=d'(d,\mathcal U)$
\begin{itemize}
\item[a)] $d'\ge d$
\item[b)] $d'|U\times U=d|U\times U$
\item[c)] If $d$ is finite--valued on every $U\in\mathcal U$ then
  $d(x,y)<\infty$ if and only if there is a stick path from $x$ to $y$.
\end{itemize}
}
\end{para}

\begin{para}\label{para8.3}
{\em An alternative way to describe this construction is the following: Let
$\Gamma_{\mathcal U}$ be the following graph. The vertices of
$\Gamma_{\mathcal U}$ are the points of $X$ and the edges of
$\Gamma_{\mathcal U}$ are the sticks, i.e., the pairs $(x,y)$
contained in one $U\in\mathcal U$. So the graph $\Gamma_{\mathcal
  U}$ is closely related to the nerve of the covering $\mathcal
U$.
To every edge $(x,y)$ of $\Gamma_{\mathcal U}$ we can associate
the weight $d(x,y)$. Then for points $x,y$ in $X$ the pseudometric
$d'(x,y)$ is the graph distance of the corresponding vertices of this
weighted graph. }
\end{para}

Let us now return to the case we are interested in. Thus, let $X$ be a
$\sigma$--compact space with a proper action of a locally compact
topological group $G$. Let $F$ be an open fundamental set for $G$ in
$X$. We consider the covering $\mathcal U$ by the translates of $F$,
so $\mathcal U=\{g F\ ;\ g\in G\}$. We apply the measuring stick
construction for an appropriate pseudometric $d$ and show that the
resulting pseudometric $d'$ is proper, but may be infinite--valued. We
do this first for the case that the orbit space $G\backslash X$ is
compact and then for the general case. We shall need an auxiliary
result about Lebesgue numbers of our covering, see below lemma
\ref{lem8.5}.The problem of infinite values of
$d'$ will be dealt with in the next section. The method will be the ``bridge
construction''.

We start with a well known result, for which we include
a proof for the convenience of the reader.

\begin{lemma}\label{lem8.4}
If the orbit space $G\backslash X$ is compact then every fundamental
set is relatively compact. Conversely, if $G\backslash X$ is compact
then every relatively compact subset $F$ of $X$ with the property that
$GF=X$ is a fundamental set for $G$ in $X$.
\end{lemma}

\begin{proof}
The second claim is clear, since property b) of a fundamental set
follows immediately from the hypothesis that the action of $G$ on $X$
is proper, see proposition and definition \ref{prop2.3} ii). To prove
the first claim choose
a compact neighborhood $U_x$ for every point $x\in X$. A finite
number of the $\pi(U_x)$, $x\in X$, cover $G\backslash X$, where
$\pi$ is the natural map $\pi:X\longrightarrow G\backslash X$,
which is known to be an open map. Let us say $G\backslash
X=\pi(U_{x_1})\cup\dots\cup \pi(U_{x_n})$, so $X=GU_{x_1}\cup\dots\cup
GU_{x_n}$. Hence $A\subset G_{U_{x_1},A} U_{x_1}\cup\dots\cup
G_{U_{x_n},A}U_{x_n}$ for every subset $A$ of $X$. For $A=F$ the
subsets $G_{U_{x_i},F}$ of $G$ are relatively compact, by property b)
of a fundamental set, see definition \ref{def6.1}. Hence $F$ is
relatively compact.
\end{proof}

A family $\mathcal D$ of pseudometrics is called {\em saturated} if
$d_1$, $d_2\in\mathcal D$ implies $\sup (d_1,d_2)\in\mathcal D$.

\begin{lemma}\label{lem8.5}
Let $\mathcal D$ be a saturated family of $G$--invariant
pseudometrics inducing the topology of $X$. Suppose the
orbit space $G\backslash X$ is compact. Then there is a
pseudometric $d\in\mathcal D$ and a positive number $\epsilon$ such
that for every $x\in X$ the ball $B_d(x,\epsilon)$ is contained in one
translate of $F$.
\end{lemma}

A number $\epsilon$ with this property is called a {\em Lebesgue
  number} for the covering $\{g F; g\in G\}$ with respect to $d$.

\begin{proof}
By $G$--invariance, it suffices to show this for points $x\in
F$. Since $\overline F$ is compact, it is covered by a finite number
of $gF$, say $\overline F\subset g_1 F\cup\dots\cup g_n F$.
Recall that $F$ is supposed to be open.
The set of
balls $B_d(x,r)$, $d\in\mathcal D$, $x\in X$, $r>0$, form a base of
the topology of $X$, not only their finite intersections, since
$\mathcal D$ is saturated. Thus there is for every $x\in\overline F$ a
pseudometric $d_x\in\mathcal D$ and a radius $r_x$ such that
$B_{d_x}(x,r_x)$ is contained in one translate of $F$, since $F$ is
open. A finite number of balls $B_{d_x}\left(x,\frac{r_x}2\right)$
  cover $\overline F$, say those for $x=x_1,\dots,x_n$. Thus for every
  $y\in\overline F$ there is an $x_i$, $i=1,\dots,n$, such that $y\in
  B_{d_{x_i}}\left (x_i,\frac{r_{x_i}}2\right)$ and hence
    $B_{d_{x_i}}\left (y,\frac{r_{x_i}}2\right )\subset
    B_{d_{x_i}}(x_i,r_{x_i})$ is contained in one translate of
    $F$. Hence our claim holds for
    $d=\sup(d_{x_1},\dots,d_{x_n})\in\mathcal D$ and $\epsilon =\inf
    (r_{x_1},\dots,r_{x_n})$.
\end{proof}

Now let again $\mathcal U=\{g F\ ;\ g\in G\}$ and for a
$G$--invariant pseudometric $d$ on $X$ let $d'=d'(d,\,\mathcal U)$ be
the pseudometric obtained by the measuring stick construction.

\begin{proposition}\label{prop8.6}
Suppose the orbit space $G\backslash X$ is compact. Let $d$ be a
continuous $G$--invariant pseudometric on $X$, for which
there is a Lebesgue number for $\mathcal U$. Then $d'$ is a proper
pseudometric, i.e., $B_{d'}(x,R)$ is relatively compact for every $x\in
X$ and every $R<\infty$.
\end{proposition}

\begin{proof}
We may assume that $x\in F$, by $G$--invariance. Then $y\in
B_{d'}(x,R)$ if and only if there is a stick path $x=x_0,
x_1,\dots,x_n=y$ with $d$--length $\sum^n_{i=1}d(x_{i-1},x_i)<R$. We
may assume that no three consecutive points $x_{i-1}$, $x_i$, $x_{i+1}$ of
our stick path are contained in one translate of $F$, because
otherwise we can leave out $x_i$ from our stick path and obtain a
stick path of not greater $d$--length. Let $\epsilon$ be the Lebesgue
number for $\mathcal U$ with respect to $d$. It follows that
$d(x_{i-1},x_i)+d(x_i,x_{i+1})\ge\epsilon$ for every $i=1,\dots,n-1$,
because otherwise $x_{i-1}$, $x_i$, $x_{i+1}$ are contained in one
translate of $F$. We thus obtain the following upper bound for the
length $n$ of our stick path:
$$
n<\frac{2R}\varepsilon + 1.
$$
Thus, let $N\in\mathbb N\cup\{0\}$ and let $B_N$ be the set of points $y\in X$
for which there is a stick path of length $N$ starting at a point
$x\in F$ and ending at $y$. We have to show that $B_N$ is relatively
compact for every $N\in\mathbb N\cup\{0\}$. For $N=0$ we have $B_N=F$. If $y\in
B_{N+1}$ there is a point $y'\in B_N$ such that $(y',y)$ is a stick,
say $\{y',y\}\subset g\,F$. Then $y'\in B_N\cap g\,F$ and hence $g\in
G_{F,B_N}=G_{B_N^{-1}, F}$. This subset of $G$ is relatively compact
by induction and property b) of a fundamental set. Thus $y\in
g\,F\subset G_{F,B_N}F$, hence $B_{N+1}\subset G_{F,B_N} F$ and so
$B_{N+1}$ is relatively compact.
\end{proof}

This yields step 4 of our proof for the case that the orbit space is
compact. For the general case we need {\em one} pseudometric $d$ for
which there is a Lebesgue number for every subset of $X$ of the form
$\pi^{-1}(K)$ where $K$ is a compact subset of $G\backslash X$. Here
we have to suppose that the orbit space is $\sigma$--compact.

Before we proceed to do this we need to figure out where $d'$ is
finite. Let $F$ and $\mathcal U$ be as above. We do not suppose that
the orbit space is compact. Let $d$ be a $G$--invariant pseudometric
on $X$ for which $d|F\times F$ has finite values. Let the symbol
``$\sim$'' denote the smallest $G$--invariant equivalence relation on
$X$ for which $F$ is contained in one equivalence class. Recall that
$G_{FF}=\{g\in G;\ g F\cap F\ne\emptyset\}$. Let $G_0$ be the subgroup
of $G$ generated by $G_{FF}$.

\begin{lemma}\label{lem8.7}
Let $x$ and $y$ be points of $X$. The following properties of the pair
$(x,y)$ are equivalent
\begin{itemize}
\item[a)] $d'(x,y)<\infty$.
\item[b)] There is a stick path from $x$ to $y$.
\item[c)] $x\sim y$.
\item[d)] The vertices $x$ and $y$ of the graph $\Gamma_{\mathcal U}$
  belong to the same connected component of $\Gamma_{\mathcal U}$.
\item[e)] If $x\in g\,F$ and $y\in h\,F$ then $g^{-1} h\in G_0$.
\end{itemize}
\end{lemma}

The equivalence classes will be called {\em islands} from now on.

\begin{proof}
a)$\Longleftrightarrow$b) was noted above, and
b)$\Longleftrightarrow$d) and b)$\Longleftrightarrow$c) follow
immediately from the definitions. \\
b) $\Longrightarrow$ e). Let $x\in
g\,F$ and $y\in h\,F$ and let $(x,y)$ be a stick, say $\{x,y\}\subset
k\,F$ for some $k\in G$. Then $g^{-1}k\in G_{FF}$ and $h^{-1}k\in
G_{FF}$ hence $g^{-1}h\in G_0$. The claim b)$\Longleftrightarrow$e)
follows now by induction on the length of the stick path.

e) $\Longrightarrow$ c). Let $Y$ be an equivalence class of
$\sim$. Thus, if one point of a translate $g\,F$ of $F$ is contained
in $Y$ then $g\,F$ is contained in $Y$. By the same argument
applied to $g\,k\,F$ with $k\in G_{FF}$ it then follows that $g\,G_{FF}
F\subset Y$, hence $g\cdot G_{FF}\cdot G_{FF}\, F\subset Y$, etc. So
$g\,G_0\, F\subset Y$ if $g\,F\cap Y\ne\emptyset$, which proves
our claim.
\end{proof}

\begin{corollary}\label{cor8.8}
The map $g\,G_0\longmapsto g\,G_0\,F$ establishes a bijection between
the set $G/G_0$ of left cosets of $G_0$ in $G$ and the set of
islands in $X$.
\end{corollary}

\begin{corollary}\label{cor8.9}
If $G\backslash X$ is $\sigma$--compact, so are $\overline F$,
$G_{\overline F,\overline F}$, $G_0$ and every island.
\end{corollary}

\begin{proof}
If $K$ is a compact subset of $G\backslash X$, then so is
$F_K:=\overline F\cap\pi^{-1}(K)=\overline{\pi^{-1}(K)\cap F}$, by
lemma \ref{lem8.4}, and hence also $G_{F_K,F_K}$, since the action of $G$
on $X$ is proper and continuous. It follows that if $G\backslash X$ is
$\sigma$--compact, so are $\overline F$, $G_{\overline F,\overline
  F}$, the subgroup $G_1$ of $G$ generated by $G_{\overline
  F,\overline F}$ and $G_1\overline F$. It thus remains to
be shown that $G_0=G_1$ and $G_0\overline F=G_0 F$. But
clearly $G_{\overline F F}= G_{FF}$ since $F$ is open, hence
$G_{\overline F\overline F}\subset G^{-1}_{\overline F F}\cdot
G_{\overline FF}$, by the formula following definition \ref{def6.1},
and thus $G_{\overline F,\overline F}\subset G_0$ and hence
$G_1=G_0$. Furthermore $\overline F\subset G^{-1}_{\overline FF}F$, by
\ref{def6.1} a), and hence $G_0\overline F=G_0 F$.
\end{proof}

We come back to the Lebesgue number and show properness of $d'$ for
the case that the orbit space is $\sigma$--compact. This accomplishes
step 4 of our plan in section \ref{sec4}. Note that at this point we
do not need that $X$ is $\sigma$--compact, only that the orbit space
is $\sigma$--compact.

\begin{lemma}\label{lem8.10}
Suppose the orbit space $G\backslash X$ is $\sigma$--compact.
\begin{itemize}
\item[a)] Then there is a continuous orbitwise proper $G$--invariant
  pseudometric $d$ on $X$ with the following properties: $d$
 is finite--valued on every island and for every compact subset $K$ of
  $G\backslash X$ there is a Lebesgue number for the covering
  $\mathcal U|\pi^{-1}(K)$ of the $G$--space $\pi^{-1}(K)$ with
  respect to the restriction of $d$ to $\pi^{-1}(K)$.
\item[b)] If $d$ is as in a) then $d'$ is proper, which means that the
  ball $B_{d'}(x,R)$ has compact closure for every $x\in X$ and every
  $0<R<\infty$.
\end{itemize}
\end{lemma}

\begin{proof}
a)\ Let $K_n$, $n\in\mathbb N$, be a sequence of compact subsets of
$G\backslash X$ such that $\bigcup^\infty_{n=1}K_n=G\backslash X$ and
$K_n\subset \overset\circ K_{n+1}$ for every $n\in\mathbb N$. Put
$X_n=\pi^{-1}(K_n)$. Then $X_n$ is a closed $G$--invariant subset of
$X$ on which $G$ acts properly with compact orbit space $K_n$. The set
$F_n:=F\cap X_n$ is an open fundamental set for $G$ in $X_n$, hence
relatively compact in $X_n$ and in $X$. So there is a continuous
orbitwise proper $G$--invariant finite--valued pseudometric $d_n$ on
$X$ such that there is a Lebesgue number for the covering $\{g
F_n;g\in G\}$ of $X_n$ with respect to the pseudometric $d_n$
restricted to $X_n$. Note that $d_n$ is defined and finite--valued on
all of $X$. To see the existence of such a $d_n$, we apply lemma
\ref{lem8.5} to the family $d|X_n\times X_n$ where $d$ runs through a
saturated family of finite--valued $G$--invariant pseudometrics on $X$
inducing the topology of $X$, which we may assume to be orbitwise
proper, by Step 3 in section \ref{sec7}.

Let $Y$ be the island $G_0 F$ containing $F$. We use here the
notation of lemma \ref{lem8.7} and its corollaries. Since $Y$ is
$\sigma$--compact, there is a family $L_n$, $n\in\mathbb N$, of
compact subsets of $Y$ such that $\bigcup^\infty_{n=1}L_n=Y$ and
$L_n\subset \overset\circ L_{n+1}$. We may assume that $d_n|L_n\times L_n$ has
values $\le 1$, by rescaling. Now define
$$
d(x,y) =
\begin{cases}
\Sigma\frac1{2^n}d_n(x,y)\ & \mathrm{if }\ x\sim y\\
\infty & \mathrm{otherwise.}
\end{cases}
$$
Then $d$ is $G$--invariant continuous orbitwise proper pseudometric
on $X$, which is finite--valued on $Y\times Y$ and hence on every
island. There is a Lebesgue number for the covering
$\{g\,F_n\,;\,g\in G\}$ of $X_n$ with respect to $d$, since there is
one for $d_n$ and $d\ge\frac1{2^n}d_n$. Here we think of $d$ and $d_n$
as restricted to $X_n\times X_n$. This implies our claim under
a).

b) Islands are of the form $g\,G_0\, F$, hence
open, since $F$ is supposed to be open. It follows that they are also
closed. Again, let $Y=G_0 F$ be the island containing $F$. Let
$B_{d'}(x,R)$, $x\in X$, $0<R<\infty$, be a ball for the pseudometric
$d'$ and let $B$ be its closure. We have to show that $B$ is
compact. We know that $K:=\pi(B)$ is compact, since $d$ is orbitwise
proper and hence so is $d'$, since $d'\ge d$ by \ref{para8.2} a). We
may assume that $x\in F$ and hence $B_{d'}(x,R)\subset Y$ and thus
$B\subset Y$.

The subgroup $G_0$ of $G$ is open since generated by the open subset
$G_{FF}$. It follows that $G_0$ is a closed subgroup of $G$.  Then the
action of $G_0$ on $Y$ is proper, since the restricted action of $G_0$
on $X$ is proper and $Y$ is a closed $G_0$--invariant subset of
$X$. And $F$ is an open fundamental set for $G_0$ in $Y$. Let $Z=Y\cap
\pi^{-1}(K)$. This is a closed $G_0$--invariant subset of $Y$ and
$F_Z:=Z\cap F=F\cap\pi^{-1}(K)$ is an open fundamental set for $G_0$
in $Z$. The orbit space $G_0\backslash Z$ is compact; it can be
identified with $K$. So we can apply proposition \ref{prop8.6} to the
$G_0$--space $Z$, the pseudometric $d|Z\times Z$ and the covering
$\mathcal U_Z :=\{g F_Z\,;\,g\in G_0\}$ to obtain that the resulting stick
path pseudometric $d'':=d'(d|Z\times Z,\mathcal U_Z)$ is proper. It
remains to see that $B_{d''}(x,R)=B_{d'}(x,R)$. Clearly $d''(x,y)<R$
implies $d'(x,y)<R$, by looking at the stick paths for $\mathcal
U_Z$. Conversely, if $d'(x,y)<R$ then there is a stick path $x=x_0$,
$x_1,\dots,x_n=y$ for $\mathcal U$ with $\Sigma
d(x_{i-1},x_i)<R$. Then all the $x_i$ are in $B_{d'}(x,R)\subset Y$
and $\pi(x_i)\in K$, hence $x_i\in Z$ and every pair $x_{i-1},x_i$ is
contained in some translate $gF$ of $F$. But then $g\in G_0$, by
\ref{lem8.7} e), and so $\{g^{-1} x_{i-1}, g^{-1}x_i\}$ is contained
in $F$ and in $Z$, hence in $F_Z$. Thus our stick path is also a stick
path for $\mathcal U_Z$ in $Z$ and thus $d''(x,y)<R$.
\end{proof}

\section{Bridges}\label{sec9}

Again, let $X$ be a $\sigma$--compact space and let the locally
compact group $G$ act properly on $X$. Note that then $G$ is
$\sigma$--compact as well, since if $X$ is the union of countably many
compact subsets $K_n$ then $G$ is the union of the countably many sets
$G_{K_n,K_n}$ which are compact since the action of $G$ on $X$ is both
proper and continuous. Let us again fix an open fundamental set $F$
for $G$ in $X$. Then, using the notation of the last section, $G_0$ is
an open subgroup of $G$ and hence $G/G_0$ is a countable discrete
space. We can thus choose a finite or infinite sequence of elements
$g_n$, $n=0,1,\dots,$ such that $G$ is the union of the disjoint
cosets $g_n G_0$. We may assume that $g_0$ is the identity
element. Let $S$ be the set of indices, so $S=\mathbb N\cup\{0\}$ or
$S=\{0,1,\dots,N\}$ for some $N\in\mathbb N\cup\{0\}$. Thus
$G=\bigcup_{n\in S} g_n G_0$ and hence $X$ is the
union of the disjoint subsets $g_n G_0 F$, $n\in S$, by corollary
\ref{cor8.8}.
Recall that the sets of the form $g\,G_0 F$ are called
islands. Consequently we define
a {\em bridge} to be a 2--point subset of $X$ of the form $\{g x, gg_n
x\}$ with $g\in G$, $n\in S$, $n\ne 0$, and $x\in F$. Note that $gx$
and $gg_n x$ are always on different islands since $n\ne 0$. But the
representation of a bridge in the form above may not be unique. Now
suppose a $G$--invariant pseudometric $d$ on $X$ is given. We then
define the {\em bridge path pseudometric} $d_B$ on $X$ as the supremum
of all pseudometrics $d''$ with the following two properties.

\begin{para}\label{para9.1} a) For every island $Y$ in $X$ we have
  $d''|Y\times Y\le d|Y\times Y$. \\
b) $d''(gx,g g_n x)\le n$ for $g\in G$, $n\in S$ and $x\in F$.
\end{para}

There is an alternative description of $d_B$ in terms of paths. Let us
define the length of a bridge $\{y,z\}$ as the smallest number $n\in S$ such
that $\{y,z\}=\{g x, g g_n x\}$ for some $g\in G$ and $x\in F$. Thus,
the length of a bridge is always an integer $\ge 1$. Let us call a
sequence of points $x=x_0,x_1,\dots,x_n=y$ a {\em bridge path} of
length $n$ from $x$ to $y$ if any two consecutive points either lie on
a common island or form a bridge, i.e., for every $i=1,\dots,n$ there
is either an island $Y$ such that $\{x_{i-1},x_i\}\subset Y$ or
$\{x_{i-1},x_i\}$ is a bridge. Define the $d$--length of such a bridge
path as $\sum^n_{i=1} d_i$ where $d_i = d(x_{i-1},x_i)$ if
$\{x_{i-1},x_i\}$ is on one island or, if $\{x_{i-1},x_i\}$ is a bridge,
then let $d_i$ be the length of this bridge.

\begin{para}\label{para9.2}
$d_B(x,y)$ is the infimum of $d$--lengths of all bridge paths from $x$
to $y$.
\end{para}

\begin{proof}
The pseudometric $d''$ defined by the statement of \ref{para9.2} has
the properties \ref{para9.1} a) and b). Conversely, if $d''$ is a
pseudometric with the properties \ref{para9.1} a) and b), then
$d''(x,y)$ is at most equal to the $d$--length of any bridge path from
$x$ to $y$, cf. the similar proof in \ref{para8.2}.
\end{proof}

\begin{proposition}\label{prop9.3}
{\bf Properties of} $d_B$
\begin{itemize}
\item[a)] $d_B$ is $G$--invariant.
\item[b)] $d_B$ is finite--valued if $d|Y\times Y$ is finite--valued for
  one (equivalently every) island $Y$.
\item[c)] If $x$ is a point of the island $Y$, then the balls
  $B_d(x,r)\cap Y$ and $B_{d_B}(x,r)$ coincide for $r<1$.
\item[d)] If $d$ is continuous, so is $d_B$.
\item[e)] Suppose $d$ is continuous, proper and, for every island $Y$,
  has finite values on $Y\times Y$. Then $d_B$ is continuous, proper
  and finite--valued (everywhere).
\end{itemize}
\end{proposition}

\begin{proof}
\begin{itemize}
\item[a)] follows from our construction.
\item[b)] follows from the fact that $d_B$ is $G$--invariant and every
  island can be reached from $F$ by a bridge.
\item[c)] follows from \ref{para9.2} and the fact that every bridge
  has length $\ge 1$.
\item[d)] A pseudometric is continuous if it is continuous near the
  diagonal, by the triangle inequality. So d) follows from c).
\item[e)] is the main point of these properties. It remains to be
  shown that $d_B$ is proper if $d$ is proper, continuous and on
  every island finite--valued. Thus let $x\in X$ and $0<R<\infty$. We
  have to show that
  $B_{d_B}(x,R)$ has compact closure. For a point $y\in X$ we have
  $d_B(x,y)<R$ if there is a bridge path $x=x_0,\dots,x_n=y$ with
  $d$--length $\Sigma d_i<R$. We may assume that three consecutive
  points $x_{i-1}, x_i,x_{i+1}$ of our bridge path are not on a common
  island, since otherwise we could leave out $x_i$ without increasing
  the $d$--length of our path, by the triangle inequality for $d$. So
  our path has at least $\frac{n+1}2$ bridges, all of length $\ge
  1$. We thus have an upper bound for the length $n$ of our bridge
  path, namely $n\le 2 R+1$. Furthermore, every bridge in our path has
  length at most $R$ and every step $d_i=d(x_{i-1},x_i)$ on one island
  has length at most $R$. It thus suffices to prove the following two
  claims.
\begin{itemize}
\item[a)] If $K$ is a compact subset of $X$, then $B_d(K,R)=\{y\in
  X\,;\,d(x,y)<R\}$ has compact closure.
\item[b)] If $K$ is a compact subset of $X$, then the set
  $B(K,R):=\{z\in X$; there is a bridge $\{y,z\}$ from a point $y\in
  K$ to $z$ of length $\le R\}$ has compact closure.
\end{itemize}
\end{itemize}

{\it Proof of a).}\ \ $K$ is contained in a finite union of islands,
  since $K$ is compact and the islands are open and disjoint and form a cover
  of $X$. It thus suffices to prove our claim for the case that $K$ is
  contained in one island, say $Y$. Let $x$ be a point of $K$. Then
  the function $y\longmapsto d(x,y)$ is continuous and finite--valued
  on $Y$, hence has a finite maximum on $K$, so $K\subset B_d(x,r)$
  for some $0<r<\infty$. Then $B_d(K,R)\subset B_d(x,r+R)$, which has
  compact closure by hypothesis. This shows our claim.\\[3mm]
{\it Proof of b).}\ \ The bridges $\{y,z\}$ starting from a point of
$K$ and having length $\le R$ are of the form $\{g x, gg_n x\}$ with
$x\in F$ and $n\le R$, and either $gx\in K$ or $gg_n x\in K$. Hence
$g\in G_{FK}$ or $g\in G_{g_n F,K}=G_{FK}\cdot g^{-1}_n$ and hence the
endpoint $z$ of our bridge is of the form $z=gg_n x\in G_{FK} g_n K$
in the first case or of the form $z=gx\in G_{FK} g_n^{-1}K$ in the
second case, thus every endpoint $z$ of such a bridge is contained in
the relatively compact set $\bigcup_{n\le R}G_{FK} g_n^{\pm 1}K$, as
was to be shown.
\end{proof}

\begin{para}\label{para9.4}
{\rm We are now ready to finish the proof of our main theorems
\ref{theo1.1} and \ref{theo1.2}. Let $X$ be a $\sigma$--compact
Hausdorff space
and suppose the locally compact topological group $G$ acts properly on
$X$. We have shown that then there is a family of continuous
$G$--invariant pseudometrics on $X$ inducing the topology of $X$, see
step 2 in chapter \ref{sec6}, which we may furthermore assume to be
finite--valued and orbitwise proper, by step 3 in chapter \ref{sec7}. Then the
stick construction of chapter \ref{sec8} gave us a pseudometric, which
is continuous, proper and on every island finite--valued, namely the
pseudometric $d'$ of lemma \ref{lem8.10}. Continuity of $d'$ follows
from property \ref{para8.2} b) and finiteness on islands from lemma
\ref{lem8.7}. If we use this pseudometric in the bridge construction
of chapter \ref{sec9} then the resulting pseudometric $d_B$ is
continuous, finite--valued and proper. If now $\mathcal D$ is a family
of $G$--invariant pseudometrics inducing the topology of $X$ -- we know
that such a family exists, by step 2 in chapter \ref{sec6} -- then the
family $\{\sup(d,d_B)\,;\,d\in\mathcal D\}$ has all the properties we
want in theorem \ref{theo1.2} (theorem \ref{theo4.1}). If $X$ is
furthermore metrizable, then there is a compatible $G$--invariant
metric $d$ on $X$, by step 2 in chapter \ref{sec6}. Again, there is a
pseudometric $d_B$ which is continuous, proper, finite--valued and
$G$--invariant. Then the metric $\sup(d,d_B)$ has all these
properties, too, and is furthermore a compatible metric. This proves
theorem \ref{theo1.1} (theorem \ref{theo4.2}).
}
\end{para}

Let us point out the following corollary, due to Haagerup and Przybyszewska \cite{haag}.

\begin{corollary}\label{cor9.5}
Every second countable locally compact group has a left invariant
compatible proper metric.
\end{corollary}

\begin{proof}
The underlying space of such a group $G$ is metrizable and
$\sigma$--compact, by corollary \ref{cor5.4}. The action of $G$ on
itself by left translations is obviously proper, so there is a
compatible left invariant proper metric on $G$, by theorem \ref{theo1.1}.
\end{proof}

As a special case we obtain the following old result of Busemann
\cite{buse}.

\begin{corollary}\label{cor9.6}
The group of isometries of a proper metric space admits a compatible left invariant proper metric.
\end{corollary}

\begin{proof}
The group $G$ of isometries of a proper metric space is locally
compact and Hausdorff, see theorem \ref{theo3.2}, and second
countable, see \cite[Ch. X, \S 3.3 Corollary]{bour}, which implies our
claim by the previous corollary.

\end{proof}

\section{Concluding remarks}

In this chapter we discuss applications and related work, mention open
questions and make other remarks.

\begin{para}\label{para10.1}
{\rm In the non--equivariant context, i.e., if we consider just the
topological space $X$ without any group action, it is well known that
a $\sigma$--compact locally compact metrizable space has a compatible
proper metric, see corollary \ref{cor5.4}. More precisely, in
\cite{w-j} it is proved that if $d$ is a complete metric on such a
space $X$ then there is a proper metric on $X$ which is locally
identical with $d$, i.e., for every point $x\in X$ there is a
neighborhood of $x$ where the two metrics coincide. Note that in our
construction the metric is not changed locally in steps 4 and 5 of
chapter \ref{sec4}. Thus in the situation of theorem \ref{theo1.1} if
$d$ is a compatible $G$--invariant metric on $X$ which is orbitwise
proper then there is a $G$--invariant compatible proper metric on $X$
which is locally identical with $d$. One may thus ask the following
question: Suppose, in the situation of theorem \ref{theo1.1}, we are
given a $G$--invariant complete compatible metric on $X$. Is there a
$G$--invariant proper (compatible) metric on $X$ which is locally
identical with $d$?
}
\end{para}

\begin{para}\label{para10.2}
{\rm Given an isometric action of a group $G$ on a $\sigma$--compact
locally compact metric space $X$ with metric $d$, it is not true in
general that there is a compatible proper metric $d_p$ for which the
action of $G$ is isometric. For an example let $X=\{(x,y)\in\real^2;
  x=0\mbox{ or } x=1\}$ endowed with the metric $d=\min\{d_E,1\}$
  where $d_E$ is the Euclidean metric of $\real^2$ restricted to
  $X$. Let $G$ be the group of isometries of $(X,d)$. There is no
  compatible proper metric $d_p$ on $X$ for which $G$ acts
  isometrically, for the following reason. The group $H$ of isometries
  of $(X,d_p)$, endowed with the compact open topology, acts properly,
  hence the isotropy group $H_{(0,0)}$ of the point $(0,0)$ is compact
    and hence has compact orbits. On the other hand, let $G_{(0,0)}$
    be the isotropy group of the point $(0,0)$ in $G$. The orbit
    $G_{(0,0)}(1,0)$ of $(1,0)$ is $\{1\}\times\real$ and is not relatively
    compact in $X$. So $G$ is not contained in $H$. The point of the
    example is that the action of $G$ is not proper, no matter which
    topology we put on $G$.
}\end{para}

\begin{para}\label{para10.3}
{\rm Let us consider the following question. Under which conditions is it
true that given a compatible metric $d$ on a locally compact
$\sigma$--compact space $X$ there is a compatible proper metric $d_p$
with the same group of isometries? A sufficient condition was given by
Janos \cite{jan}, namely if $(X,d)$ is a connected uniformly locally
compact metric space.
}
\end{para}

\begin{para}\label{para10.4}
{\rm Note that if we have a closed subgroup $G$ of the group of
  isometries of a proper metric space $(X,d)$ then it is not true in
  general that there is a metric $d_1$ on $X$ for which $G$ is the
  precise group of isometries. E.g., the space $X=\real$ of real
  numbers with the Euclidean metric has the group $G=\real$ as a
  closed subgroup of its group of isometries. But for every
  $G$--invariant metric $d_1$ on $X$ we have $d_1(x,0)=d_1(0,-x)$,
  hence the group of isometries of $d_1$ contains the reflections of
  $\real$ and is thus strictly larger than $\real$.
}
\end{para}

\begin{para}\label{para10.5}
{\rm Given a proper action of a locally compact topological group $G$
  on a locally compact metrizable space $X$, one can ask if there is a
  $G$--invariant metric. This is known to be equivalent to
  $G\backslash X$ being
  paracompact \cite{koszul}, \cite{AS}, \cite{AN}. The answer is
  positive in many cases, see \cite{AS}, \cite{AN}. If $X$ is no
  longer locally compact, the answer is known to be negative if the
  action is Bourbaki--proper, see \cite{AS}, but again unknown in
  general for Palais--proper actions.
}
\end{para}

\begin{para}\label{para10.6}
{\rm Our theorem \ref{theo1.1} has potential applications for the
  Novikov conjecture. Namely, let $G$ be a locally compact
  second countable group and let $\mu$ be a Haar measure on $G$. Then,
  using a proper left invariant compatible metric on $G$, Haagerup and
  Przybyszewska have proved in \cite{haag} that there is a proper
  affine isometric action of $G$ on some separable strictly convex
  reflexive Banach space. Kasparov and Yu have recently proved that the
  Novikov conjecture holds for every discrete countable group which
  has a uniform embedding into a uniformly convex Banach space, see \cite{k-y}
}
\end{para}

\newpage

%\vglue1cm
{\sc Herbert Abels}\\[2mm]
Universit\"at Bielefeld\\
Fakult\"at f\"ur Mathematik\\
Postfach 100\,131\\
33501 Bielefeld\\
Germany\\
e--mail:
abels$@$math.uni--bielefeld.de\\[5mm]
{\sc Antonios Manoussos}\\[2mm]
Universit\"at Bielefeld\\
Fakult\"at f\"ur Mathematik, SFB 701\\
Postfach 100\,131\\
33501 Bielefeld\\
Germany\\
e--mail:
amanouss@math.uni-bielefeld.de\\[5mm]
{\sc Gennady Noskov}\\[2mm]
Sobolev Institute of Mathematics\\
Pevtsova 13\\
Omsk 644099\\
Russia\\
and\\
Universit\"at Bielefeld\\
Fakult\"at f\"ur Mathematik, SFB 701\\
Postfach 100\,131\\
33501 Bielefeld\\
Germany\\
e--mail:
g.noskov@gmail.com

\end{document}